\documentclass[12pt]{article}
 \usepackage[dvipsnames]{xcolor}
\usepackage[utf8]{inputenc}
\usepackage{amsmath,amssymb,amsthm}
\usepackage{hyperref}
\usepackage{graphicx}
\def \nn  {\nonumber}
\def\bb{{\rm\bf \textcolor{blue}{B}}}
\def\ff{{\rm\bf \textcolor{red}{F}}}
\def\gg{{\rm\bf \textcolor{ForestGreen}{G}}}

\newtheorem{coro}{Corollary}
\newtheorem{rema}{Remark}
\newtheorem{defi}{Definition}
\newtheorem{lemma}{Lemma}
\newtheorem{prop}{Proposition}
\newtheorem{thm}{Theorem}
\newcommand{\dt}{\,\mathrm{d}t}
\def \a   {\alpha}

\def \g   {\gamma}
\def \d   {\delta}
\def \l   {{\lambda}}
\def \om   {{\omega}}
\def \eps {\varepsilon}
\def\E{{\mathbb{E}}}
\def\P{{\mathbb{P}}}
\def\R {{\mathbb{R}}}
\newcommand{\Z}{\mathbb{Z}}

\def\|{\,|\,}

\def\bn#1\en{\begin{align*}#1\end{align*}}
\def\bnn#1\enn{\begin{align}#1\end{align}}

\title{Generalizations of forest fires with ignition at the origin}
\author{Francis Comets\footnote{Universit\'e de Paris, LPSM, F-75006 Paris, France},\ \ Mikhail Menshikov\footnote{Department of Mathematical Sciences, Durham University, DH1 3LE, England},\ \ Stanislav Volkov\footnote{Centre for Mathematical Sciences, Lund University, SE-22100, Sweden}} 

\usepackage{marvosym}
\usepackage{pifont}

\begin {document}
\maketitle
\begin {abstract}
We study generalizations of the Forest Fire model introduced in~\cite{BTRB} and~\cite{Volk} by allowing the rates at which the trees grow to depend on their location, introducing long-range burning, as well as a continuous-space generalization of the model. We establish that in all the models in consideration the expected time required to reach site at distance $x$ from the origin is of order $(\log x)^{(\log  2)^{-1}+\delta}$ for any $\delta>0$.
\end {abstract}



\section{Introduction}
The purpose of this paper is to generalize the results for the version of forest fire model studied in~\cite{Volk}; see also~\cite{BTRB,CRANE}. In that model, one considers the following continuous time process  on $\Z_+=\{0,1,2,\dots\}$. Let $\eta_x(t)\in \{0,1\}$ be the state of site $x\in \Z_+$ at time $t\ge 0$. Site $x$ is declared {\it vacant}  ({\it occupied} resp.) if $\eta_x=0$ ($\eta_x=1$ resp.) The process evolves as follows: vacant sites become occupied independently with rate $1$; after they are occupied, they can be ``burnt'' by a fire spread from a neighbour on the left, which makes them vacant again. There is a constant (and {\em the only}) source of fire attached to site $0$; so when site $0$ becomes occupied, the whole connected cluster of occupied sites which contains~$0$ is instantaneously burnt out. 

Initially all the sites are vacant. We are interested in the dynamics of process $\{\eta_x(t)\}$ as $t\to\infty$. For other relevant models on forest fire we refer the reader to~\cite{VDBa,VDBb,BF}; also, some more recent results dealing with complete graphs can be found in~\cite{CRANE} and for planar lattices in~\cite{KISS}. See also~\cite{MARTIN} for the connection between multiplicative coalescent with linear deletion and forest fires.

In the current paper we consider the following two generalizations.

\begin{itemize}
\item[(i)]
We allow the rates at which vacant site $x$ becomes occupied to depend on $x$, though they all must lie in a certain interval, and also  consider a long-range mode (the terminology is consistent with use in the percolation theory), where the fires can spread further than just to the immediate neighbour of the ``burning'' site. 

Note that for each site $x$ the sequence of burning times at $x$ is a renewal process that is measurable with respect to the filtration generated by the arrival processes at the sites $\{0,1,\dots,x\}$.

The results are presented in Theorem~\ref{thm_main}  below.

\item[(ii)]
We also consider a continuous-space generalization of the process where we replace $\Z_+$ by $\R_+$; the main result here is our Theorem~\ref{thm_cont}. We want to mention that in this case we consider only a homogeneous version, as even then the arguments become quite complicated. 
\end{itemize}

We note that the bounds on the expected time to reach point $x$ given by Theorems~\ref{thm_main} and~\ref{thm_cont} are probably suboptimal; we expect them to be in reality of order $\log x$ (as was shown for a  simpler model in~\cite{Volk}) rather than some power of $\log x$, however, we do not see how to prove it in a general case.

The rest of the paper is organized as follows.
In Section~\ref{sec-formal} we present the formal definition of the processes we study. 
In Section~\ref{sec-nonh} we formulate and prove our main theorem for model (i). The strategy of the proof is roughly the following. First, we introduce an auxiliary ``green'' process, which is a version of the forest process {\em without fires}, and we study how far a {\em potential} fire could reach, should it occur at the specific point of time. It turns out that it is much easier to study this auxiliary process, and at the same time the ``real'' fire process and the green process coincide from time to time. This allows us to introduce an increasing sequence of locations $n_k\uparrow \infty$, and using intricate coupling to get an upper bound of the expected time when $n_k$s are burned for the first time.
Finally, in Section~\ref{sec-cont} we formulate and prove the theorem about the continuous version of the model.

Throughout the paper we will use notation $a(t)\asymp b(t)$ if $\lim_{t\to\infty} a(t)/b(t)=1$, and we assume that all the processes are c\`adl\`ag, that is, at the time of the fire all points of the burning cluster become vacant.

\section{Formal definitions and the ``green'' process}\label{sec-formal}
The probability space on which we define the processes is the following.
Given a deterministic sequence $\{\l_x\}_{x\in\Z_+}$ for each $x$ we have an independent Poisson process of rate $\l_x$, denoted by ${\cal P}_x(t)$, started at zero. The probability and expectations throughout the paper will be with respect to the product measure generated by these processes;  elementary outcomes will be denoted by $\om$.
Note that this definition will be slightly modified for model (ii).

\subsection*{The green process for the forest fire process}
Consider the following modification of the forest fire process.

Fix an $x\in\{1,2,\dots\}$. Assume that there are no fires at all, and wait for the first time when~$x$ is {\em reachable}  from $0$  (precisely defined immediately after~\eqref{Sdef}). This will provide a trivial lower bound for the minimum time necessary for the fire to reach this point. 

More formally, let ${\cal P}_x$, $x\in \Z_+$, be a collection of independent Poisson processes, such that ${\cal P}_x$ has  rate $\l_x$. Define a $\{0,1\}^{\Z_+}$--valued continuous time stochastic process $S^\gg$ as follows.
The value $S^\gg(x,t)\in\{0,1\}$ for $x\in \Z+$, $t\ge 0$, is the state of site $x$ at time $t$. We say that $x$ is occupied at time $t$ (vacant resp.) whenever $S^\gg(x,t)=1$ ($S^\gg(x,t)=0$ resp.). At time $t=0$ all the states are vacant, i.e., $S^\gg(x,0)=0$ for all  $x\ge 0$. 

At each arrival time of the Poisson process ${\cal P}_x$, the state $x$ becomes $1$, regardless of its previous value. Formally,
\begin{align}\label{Sdef}
S^\gg(x,t)=\begin{cases}
0, &\text{if } {\cal P}_x(t)=0;\\
1, &\text{if } {\cal P}_x(t)\ge 1.
\end{cases}
\end{align}
Let $r\ge 1$ (called {\em the range} of the process) be some positive integer; we say that $x\in\{1,2,\dots\}$ is {\em    reachable from $0$ at time $t$}, if there exists a positive integer $n$ and a sequence of points in $\Z_+$  $0=x_0<x_1<x_2<\dots<x_n= x$  such that  $x_i-x_{i-1}\le r$ for $i=1,2,\dots,n$ and also $S^\gg(x_i,t)=1$ for $i=0,1,2,\dots,n-1$ (note that we {\em do not require} that $S^\gg(x,t)=1$).
This defines the {\em green process}. For the green process, we can define the quantity $N^\gg(t)$, which is the right-most reachable vertex at time~$t$; in particular, if $N^\gg(t)=n$ then $S^\gg(n-r,t)=1$ and $S^\gg(n-i,t)=0$ for $i=0,1,\dots,r-1$. 
For definiteness, if $S^\gg(0,t)=0$ let  $N^\gg(t)=0$.

The time when point $x$ becomes reachable for the first time is denoted by $\tau^\gg_x=\inf\{t>0:\ N^\gg(t)\ge x\}$.
It is easy to see that $N^\gg(t)<\infty$ a.s.\ and that $N^\gg(t)\uparrow\infty$ as $t\to\infty$ a.s.

The forest fire process, denoted by $S^\ff(x,t)$, is easily defined on the same probability space generated by the same processes ${\cal P}_x(t)$, $x=0,1,2,\dots$, as the green process, with the additional rule which says that whenever site $0$ becomes occupied, all the sites which are reachable from $0$ become vacant. We say that all these sites {\em are burnt} at this time (regardless whether they were previously occupied or not). Throughout the rest of the paper by $N(t)$ we shall denote the rightmost point burnt by the fire by time~$t$.

The previously studied case (e.g.~in \cite{Volk}) when $r=1$ is referred to as the ``short" range as opposed to the ``long" range when $r$ can be greater than $1$.

\begin{rema}
Suppose that $\l_x\equiv1$ and $r=1$ as in~\cite{Volk}. The process $N^\gg$ is then a time inhomogeneous pure jump Markov process on $\Z_+$ with unit jump rate and with a  {\rm Geom}$(e^{-t})$ jump distribution, so $N^\gg(t)\,e^{-t}$ converges in law to an exponential variable, and
$$
\P\big( N^\gg(t+\dt) =N^\gg(t)+k \ \big  \vert\  {\rm history\ up\ to\ } t\big) = 
\left\{
\begin{array}{ll}
1-\dt+o(\dt), & k=0,    \\
 e^{-t}(1-e^{-t})^{k-1} \dt + o(\dt),&  k \geq 1.   
 \end{array}
\right.
$$
Indeed, for an individual site the probability to be occupied in time $t$ is $1-e^{-t}$, and thus the longest stretch of such occupied sites has a geometric distribution. Since, given that at time $t$  we have $N^\gg(t)=n$, the site $n+1$ must be vacant at time $t$. If it becomes occupied during time $dt$, which occurs  with probability $dt$, the process will immediately spread further to the right due to the already occupied sites at $n+1,n+2,\dots$.
\end{rema}

The usefulness of the green process will become clear from Proposition~\ref{prop_green}, introduced a bit later.
At this point we outline our construction. The green and  fire processes are naturally coupled, they have the same ``reach'' at a sequence of times tending to infinity. The green process is a trivial upper bound for the other one,  while it can be analysed somewhat explicitly.  On the other hand, the fire process has enough of a renewal structure to estimate how much slower than the green process it can be.
\\[3mm]
Now we rigorously introduce models (i) and (ii) studied in the current paper.

\subsection*{(i) Non-homogeneous, long-range process}
Suppose that the rate at which site $x$ become occupied depends on the site, and denoted by $\l_x$. Throughout the paper we will assume uniform bounds on these rates
$$
c_1\le \lambda_x \le c_2,\qquad  x=0,1,2\dots
$$
for some fixed constants $c_2>c_1>0$.  We will assume also that the process is long range, that is, when the fire hits point $x$ it burns all vertices in $[x,x+1,\dots,x+r]$ for some fixed positive integer $r$ (``range''). 

In the original model studied in~\cite{Volk}, one has $\l_x= 1$ for all $x$, and $r=1$. We will refer to the model where all $\l_x$ are the same as to the {\em homogeneous} model; in case where $r=1$ we will say that it is a {\em short-range} model.

\subsection*{(ii) Continuous tree model on $\R_+$}

\setlength{\unitlength}{10mm}
\begin{picture}(17,2)(0,0)
\put(1,1){\line(1,0){10}}
\put(1,1){\circle*{0.1}}
\put(5,1){\circle{1}}
\put(5.6,1){\circle{1}}
\put(7,1){\circle{1}}
\put(9,1){\circle{1}}
\thicklines
\put(1,1){\circle{1}}
\put(1.8,1){\circle{1}}
\put(2.3,1){\circle{1}}
\put(3.2,1){\circle{1}}
\end{picture}

\noindent
Suppose that ``trees'' (or, rather, their centres) arrive on $\R_+$ as a Poisson process of rate $1$. Namely, the number of trees which appeared on $[0,x]$, $x>0$ by time $t>0$ is given by a number of points of two-dimensional Poisson process in a rectangle $[0,x]\times [0,t]$. Each tree is a closed interval (one-dimensional circle) of a fixed radius $1$. There is a constant source of fire attached to the origin, point~$0$. Whenever a tree covering the origin appears, it immediately burns down together with the whole connected component of trees (overlapping circles) containing point~$0$. For the continuous model, we can similarly define $N(t)$ and $N^\gg(t)$, as we did for the model on $\Z_+$.

Also, to avoid making the paper unnecessarily cumbersome, we restrict our attention only the homogeneous version of model (ii), even though one can think of a more general setup.

\begin{defi}\label{def1}
Let $\tau_x$ , $x\in \R_+$, be the time when point $x$ is burnt by fire for the first time in the forest fire model (either (i) or (ii)). Similarly, let $\tau_x^\gg$ be the first time when $x$ is reachable by the green process. More formally,
\begin{align*}
\tau_x&=\inf\{t:\ N(t)\ge x\},\\
\tau_x^\gg&=\inf\{t:\ N^\gg(t)\ge x\}.
\end{align*}
\end{defi}

\begin{figure}
\begin{center}
\includegraphics[scale=1]{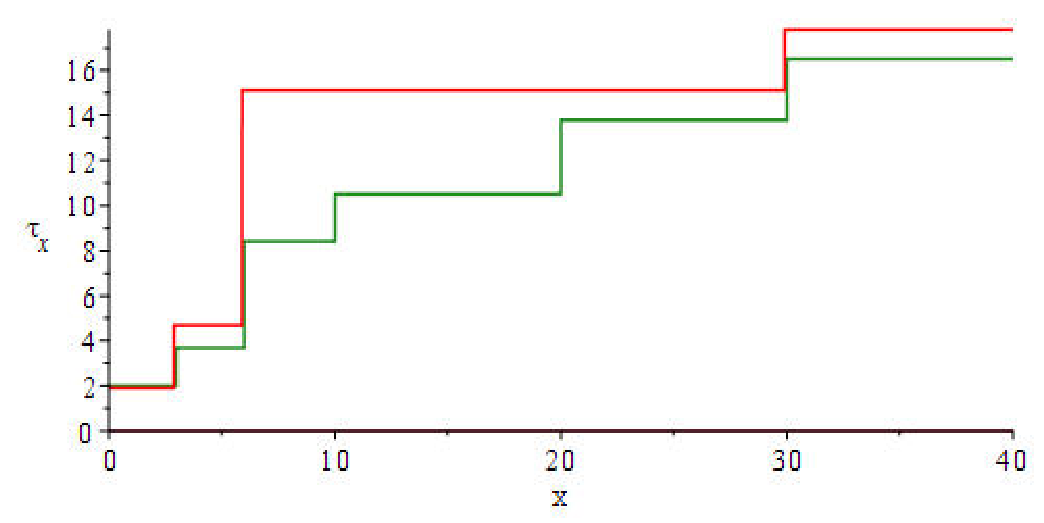}
\caption{$\tau_x$ for the forest fire and green processes}
\label{figr}
\end{center}
\end{figure}

The following statement applies to both models (i) and (ii).
\begin{prop}\label{prop_green}
Let $\tau_x$ and $\tau_x^\gg$ be as defined above. Then there is a coupling between the green process and the forest fire process such that
\begin{itemize}
\item the original forest fire process is always behind or equal to the green process, i.e.\ $\tau_x\ge \tau_x^\gg$ for all~$\omega$;
 
\item nevertheless, every time the fire process burns a point never burnt before, its ``spread'' coincides with that of the green process (see Figure~\protect\ref{figr}): namely, for almost all $\omega$ there is an infinite increasing sequence of times $t_i\to\infty$, 
$i = 1, 2,\dots$, such that $N(t_i)=N^\gg(t_i)$.
\end{itemize}
\end{prop}
\begin{proof}
The first part of the statement follows immediately from the construction of the two process on the same probability space. To show the second part, introduce the increasing sequence of stopping times $\sigma_i$ and locations $u_i$ such that $\sigma_0=u_0=0$;
$\sigma_i>\sigma_{i-1}$ is the first time when the fire burns a point never burnt before, and $u_i$ is the right-most point burnt by the fire at time $\sigma_i$. It is not hard to see that all $\sigma_i<\infty$ and $u_i<\infty $ a.s.

Then at times $\sigma_i$ the forest fire process coincides with the ``green" process, since at the infinitesimal moment before that all the vertices to the right of $u_{i-1}$ are in the same state for both processes.
\end{proof}

\section{Non-homogeneous and long-range processes}\label{sec-nonh}
The main result of this section is
\begin{thm}\label{thm_main}
Assume that the forest fire process is the one described by model (i) for some $r$, $c_1$, and~$c_2$. Fix any $\kappa>(\log 2)^{-1}=1.442695\dots$. Then, for all large enough $x$,
$$
\E \tau_x\le  (\log x)^{\kappa}.
$$ 
\end{thm}
Observe that this result does not depend on values $c_1,c_2,r$. The proof of the above statement will follow from Sections~\ref{Sec_Nonhom_long} and~\ref{Sec_couplingfire}.

\begin{rema}
By comparing  $\tau_{x+1}$  with $\tau_x$, it is not difficult to see that $\E \tau_{x+1}\le \left(2+o(1)\right)\,\E\tau_x$, and hence to get by induction that $\E\tau_x<\infty$ for all $x\ge 1$. See the proof of Theorem~\ref{thm_main}.
\end{rema}

\begin{rema}
Using Propositions~\ref{prop_green} and~\ref{propbounds} (or Propositions~\ref{prop_green} and~\ref{prop_homog} in the homogeneous case), which come later in the paper, it is easy to obtain that there is a sufficiently small constant $C>0$ such that $\P(\tau_x\ge C \log x)\ge \P(\tau_x^\gg \ge C \log x)\to 1$ as $x\to\infty$, and hence $\tau_x$ must be at least of the order $\log x$.
\end{rema}

\begin{rema}
We conjecture that the correct power of the $\log$ in the statement of Theorem~\ref{thm_main} should be, in fact, $1$, as in~\cite{Volk}, thus matching the lower bound.
Please also see Corollary~\ref{corliminf}.
\end{rema}

Recall that the ``green process'' is  the modified version of the forest fire process in which there are no fires at all, and that $N^\gg(t)$ is the size of the connected cluster of occupied vertices including the origin for the green process.  We start with a special case.

\subsection{Short range non-homogeneous process}
This is the special case of the process when $r=1$, which we can study  without reference to the green process, using Lemma~3 in~\cite{BTRB} about the  invariance of the distribution of $\tau_x$ w.r.t.\ permutation of the rates $\lambda_1,\dots,\lambda_x$.

It turns out that larger $\lambda$'s correspond to smaller $\tau_x$, hence $\tau_x$ is stochastically bounded below and above by~$\tilde\tau_x/c_2$ and $\tilde\tau_x/c_1$ respectively, where $\tilde\tau_x$ is the distribution found in~\cite{Volk}.
In particular, $\E \tau_x=O(\log x)$.

To see why this is true, note that by basic coupling the fire process with rates $\lambda_1,\ldots, \lambda_{x-1}, c_2$ hits point $x$ before $\tau_x$. Now, the permutation invariance result of~\cite[Lemma 3]{BTRB}  states that the former has same law as the hitting time of $x$ by the fire process with rates $c_2, \lambda_1,\ldots, \lambda_{x-1}$. Iterating the argument, we prove our claim.

\begin{rema}
It is straightforward that for the green process in this case we have
$$
\P(N^\gg(t)\ge n)=\left(1-e^{-\l_1 t}\right)
\left(1-e^{-\l_2 t}\right)\cdots 
\left(1-e^{-\l_n t}\right).
$$
\end{rema}

\subsection{Long-range green process}
Define $B_i=B_i(t)$ as the event that node $i$ is vacant at time $t$. Then $B_i$ are independent and  $\P(B_i)=e^{-\l_i t}$.  Moreover,
\begin{align}\label{eqNB}
\{N^\gg(t)<n\}=\bigcup_{i=0}^{n-r} \left(B_{i+1} \cap B_{i+2}\cap \dots \cap B_{i+r}\right)
\end{align}
i.e., there is a vacant interval of length $r$ somewhere on $[0,n]$.
From equation~(\ref{eqNB}) we have
\begin{align}\label{eq_alpha2}
\P\left(N^\gg(t)<n\right)\le \sum_{i=0}^{n-r} 
\P\left(\bigcap_{m=1}^r B_{i+m}\right)
=\sum_{i=0}^{n-r} e^{-(\l_{i+1}+\dots+ \l_{i+r})t} =:f_{n}(t).
\end{align}

\subsubsection{The homogeneous case}\label{subsec_homlong}
Without loss of generality assume that $\l_x=\l=1$. The next statement shows that the first time at which $N^\gg(t)\ge n$ is of order $T:=\frac{\log n}r$ where $r$ is the range of the process. 
\begin{prop}\label{prop_homog}
For every $\eps>0$ we have
\begin{align*}
\P\left( N^\gg((1-\eps)T)>n\right)&=o(1),\\
\P\left( N^\gg((1+\eps)T)<n\right)&=o(1).
\end{align*}
as $n\to\infty$.
\end{prop}

Let $p_n=p_n(t):=\P(N^\gg(t)\ge n)$, that is, the probability that on the set of vertices $\{1,2\dots,n\}$ there are no ``holes'' of length of at least $r$, and site $0$ is occupied. Observe that $p_n$ satisfies the following recursion:
\begin{align*}
p_n&=\sum_{k=1}^{r} (1-\a) \a^{k-1} p_{n-k}, \ n> r;\\
p_1&=p_2=\dots=p_{r-1}=1-e^{-t},
\end{align*}
where $\a=\a_t:=e^{-t}$. 
Indeed, in order to reach $n>r$, one must have an occupied vertex somewhere between $1$ and $r$; the probability that the first such vertex is exactly at $k$ is 
$$
\underbrace{e^{-t}\cdot e^{-t}\cdot...\cdot e^{-t}}_{k-1\text{ times}}\cdot (1-e^{-t}),
$$
and then we need to reach $n$ from $k$.
Therefore, the asymptotic behaviour of $p_n$ will depend on the largest solution of the characteristic equation
\begin{align}\label{eqchareq}
\xi^r- \sum_{k=0}^{r-1} (1-\a)\, \a^{(r-1)-k} \, \xi^{k}=0.
\end{align}
In particular, when $r=2$ we  get
$
\xi_{1,2}=\frac{1-\a\pm\sqrt{D}}{2}
$
where $D=1+2\a-3\a^2$ resulting in
$$
p_n=\frac{1}{2}\left[\left(1+\frac{1+\a}{\sqrt{D}}\right)\cdot 
\left(\frac{1-\a+\sqrt{D}}{2}\right)^n+\left(1-\frac{1+\a}{\sqrt{D}}\right)\cdot \left(\frac{1-\a-\sqrt{D}}{2}\right)^n\right],\quad n=1,2,\dots
$$

While we cannot solve~(\ref{eqchareq}) explicitly except for $r=2$ and $r=3$, we can still find the critical value of $t$ for a given $n$, that is, $T$. Indeed, from~(\ref{eq_alpha2}) we have
$$
\P\left( N^\gg(T(1+\eps))<n \right)\le n e^{-rT(1+\eps)}=n^{-\eps}\to 0,
$$
and at the same time,  since the unions of the events below are independent,
\begin{align*}
\P(N^\gg(t)>n)&\le \P\left( \left(B_1^c\cup  \dots \cup B_r^c\right) \cap
\left(B_{r+1}^c \dots B_{r+r}^c\right)\cap \dots
\cap \left(B_{Mr+1}^c\dots B_{Mr+r}^c\right)\right)
\\
&=\prod_{k=0}^{M} \P\left(\left(B_{kr+1} B_{kr+2} \dots B_{kr+r}\right)^c\right)
= \left(1-e^{-rt}\right)^{\lfloor n/r\rfloor }
\le \left(1-e^{-rt}\right)^{ n/r-1 }
\end{align*}
where $M=\lfloor n/r \rfloor -1$ and 
$\lfloor x\rfloor$ denotes the integer part of $x\in\R$. Consequently,
$$
\P(N^\gg(T(1-\eps))>n)\le  \left(1-e^{-rT(1-\eps)}\right)^{n/r-1 }
=  \left(1-\frac 1{n^{1-\eps}}\right)^{n/r-1}
\le 2\,   e^{-n^{\eps}/r}\to 0.
$$
Thus Proposition~\ref{prop_homog} is proven.

\subsubsection{The non-homogeneous case}\label{Sec_Nonhom_long}
Here we no longer assume that $\l_x$'s are the same for all $x$, as we did in Section~\ref{subsec_homlong}.

The following statement is trivial and its proof is thus omitted.
\begin{prop}\label{propfn}
The function $f_n(t)$ defined in~(\ref{eq_alpha2}) satisfies the following properties:
\begin{itemize}
\item
For a fixed $n$, $f_n(t)$ is monotonically decreasing with range from $n-r+1$ to $0$ as $t$ goes from $0$ to $+\infty$.
\item
For a fixed $t$, $f_n(t)$ is monotonically increasing in $n$ and moreover $f_{n}(t)-f_{n-1}(t)\le e^{-rc_1 t}$.
\end{itemize}
\end{prop}

Let us split the sum in $f_n(t)$ into $r$ ``almost'' equal parts according to the value of the remainder when $i$ is divided by~$r$. Namely, for $j=0,1,\dots,r-1$ let
$$
f_{n,j}(t)=\sum_{k=0}^{\lfloor \frac{n-j-r}{r}\rfloor}e^{-(\l_{rk+j+1}+\l_{rk+j+2}+\dots+ \l_{rk+j+r})t}
$$
therefore
$$
f_n(t)=f_{n,0}(t)+f_{n,1}(t)+\dots+f_{n,r-1}(t).
$$

First, we show  that 
$$
1-e^{-\frac{f_n(t)}r}\le \P\left(N^\gg(t)<n\right)\le f_n(t).
$$
Indeed, the upper bound  follows from~(\ref{eq_alpha2}). To show the lower bound note that at least one of $f_{n,j}(t)$, $j=0,1,\dots,r-1$, must be larger than $f_n(t)/r$;  w.l.o.g.\ assume $f_{n,0}(t)\ge f_n(t)/r$. Then we have
\begin{align*}
\P\left(\{N^\gg(t)<n\}^c\right)&=\P\left( \bigcap_{i=0}^{n-r} \left(B_{i+1}B_{i+2} \dots B_{i+r}\right)^c \right)
\le 
\P\left( \bigcap_{k=0}^{\lfloor n/r \rfloor-1} 
\left(
B_{rk+1} B_{rk+2}\dots B_{rk+r}\right)^c 
\right)
\\
&\stackrel{\text{indep.}}{\,=\,}
\prod_{k=0}^{\lfloor n/r \rfloor-1} 
\P\left(  \left(
B_{rk+1} B_{rk+2}\dots B_{rk+r}\right)^c 
\right)
=\prod_{k=0}^{\lfloor n/r \rfloor-1}\left( 1-e^{-(\l_{rk+1}+\dots+\l_{rk+r})t} 
\right)
\\
&\le 
\exp\left(-
\sum_{k=0}^{\lfloor n/r \rfloor-1}
e^{-(\l_{rk+1}+\dots+\l_{rk+r})t} 
\right) =e^{-f_{n,0}(t)}
\le e^{-f_n(t)/r}
\end{align*}
hence $P(N^\gg(t)<n)\ge 1-e^{-f_n(t)/r}$.

Let $t_*=t_*(n,\alpha)$ be such that $f_n(t_*)=\alpha\in(0,1)$; its existence and uniqueness follows from  Proposition~\ref{propfn} and the fact that $f_n(\cdot)$ is strictly decreasing on $[0,\infty)$ and continuous. Now set $\alpha=1/2$ and define $\tilde T=t_*(n,1/2)$; later we will establish that $\tilde T=O(\log n)$,  see~\eqref{FC1}. Then
$$
0<1-e^{-\frac{1}{2r}}\le \P\left(N^\gg(\tilde T)<n\right)\le \frac 12.
$$

\begin{prop}\label{propbounds}
Fix a small $\eps>0$. Then for any $c\in(0,c_1/c_2)$ and all large $n$
\begin{align*}
\P\left(N^\gg(\tilde T-\eps \tilde T)\ge n\right)&\le e^{-n^{c\eps}}, \\
\P\left(N^\gg(\tilde T+\eps \tilde T)<n\right)&\le n^{-c\eps}.
\end{align*}
\end{prop}
\begin{proof}
Indeed, since
$$
(n-r+1)\ e^{-rc_2 \tilde T}=\sum_{i=0}^{n-r} e^{-rc_2 \tilde T}
\le 
\sum_{i=0}^{n-r} e^{-(\l_i+\dots+\l_{i+r})\tilde T}\equiv \boxed{\frac 12}
\le
\sum_{i=0}^{n-r} e^{-rc_1 \tilde T}=(n-r+1)\ e^{-rc_1 \tilde T}
$$
by summing up the LHS, the RHS, and $1/2$ of the above chain of inequalities we immediately get
\begin{align}\label{FC1}
\left(2n-2r+2\right)^{1/c_2}\le e^{r\tilde T}\le \left(2n-2r+2\right)^{1/c_1}.
\end{align}
Consequently,
$$
f_n(\tilde T+\eps \tilde T)=\sum_{i=0}^{n-r} 
\frac{e^{-(\l_{i+1}+\dots+\l_{i+r})\tilde T}}{e^{(\l_{i+1}+\dots+\l_{i+r})\eps \tilde T}}
\le 
\sum_{i=0}^{n-r} \frac{e^{-(\l_{i+1}+\dots+\l_{i+r})\tilde T}}{e^{r\tilde T \cdot c_1 \eps }}
\le \frac{1}{2\, (2n-2r+2)^{ \frac{c_1 \eps }{c_2} }}.
$$
By the same token,
\begin{align*}
f_n(\tilde T-\eps \tilde 	T)&=\sum_{i=0}^{n-r} e^{-(\l_{i+1}+\dots+\l_{i+r})\tilde T}\times e^{
(\l_{i+1}+\dots+\l_{i+r})\eps \tilde T}\\
&\ge  \sum_{i=0}^{n-r} e^{-(\l_{i+1}+\dots+\l_{i+r})\tilde T}\times e^{r\tilde T \cdot c_1 \eps } 
\ge \frac{1}{2}\, (2n-2r+2)^{ \frac{c_1 \eps }{c_2} }.
\end{align*}
As a result
$$
\P\left(N^\gg(\tilde T-\eps \tilde T)<n\right)\ge 1-e^{-r^{-1}\,f_n(\tilde T-\eps \tilde T)}\ge 1-\exp\left(-\frac{(2n-2r+2)^{ \frac{c_1 \eps }{c_2}}}{4r}
\right)
$$
and
$$
\P\left(N^\gg(\tilde T+\eps \tilde T)<n\right)\le f_n(\tilde T+\eps \tilde T)\le \frac{1}{2\, (2n-2r+2)^{ \frac{c_1 \eps }{c_2} }}
$$
which yields the statement of the theorem, since $\frac{c_1\eps}{c_2}>c\eps$.
\end{proof}

\begin{coro}\label{corliminf}
There exists $C_1>0$ such that
\begin{align}\label{eqliminf}
\liminf_{x\to\infty} \frac{\tau_x}{\log x}\le C_1\quad\text{a.s.}
\end{align}
\end{coro}
\begin{proof}
From Proposition~\ref{propbounds} we get that for some non-random $C_1>0$
$$
\P(\tau_n^{\gg}<C_1\log n)\le n^{-c\eps}\text{ for all large }n.
$$
By choosing a sequence $n_j=j^A$, $j=1,2,\dots$ for some large $A>0$ such that $Ac\eps>1$, we get by Borel-Cantelli lemma that a.s.\ for all large enough $j$ we have
$$
\tau_{n_j}^\gg \le  C_1 \log (j^A).
$$
Using monotonicity of $\tau_x^\gg$ in $x$, and the fact that for each $x$ there is a $j$ such that 
$$
(1-o(1))\, x\le j^A\le x <(j+1)^{A}\le (1+o(1))\, x,
$$
we conclude that
$$
\tau_{x}^\gg \le  (C_1+o(1))  \log x
$$
a.s.\ for all large $x$. As a result, by the second part of Proposition~\ref{prop_green}, which equivalently stated means that $\tau_x=\tau_x^\gg$ for infinitely many $x$, we obtain~\eqref{eqliminf}.
\end{proof}

\subsection{Coupling of the forest fire  and green processes}
\label{Sec_couplingfire}
Fix $\gamma\in (1,2)$, and let us define the increasing sequence of times $\g^k$, $k=1,2,\dots$. Let
\begin{align*}
n_k&=\min\{n\in\Z_+:\ f_{n}(\g^k)\ge 1/2\};\\
T_k&=t_*(n_k,1/2).
\end{align*}

\begin{prop}\label{propTnu}
Let $n_k$ and $T_k$ be as defined above. Then the following holds.
\begin{itemize}
\item $n_k=\exp\left\{r\, \tilde c_k\, \g^k\right\}$ for some $\tilde c_k\in [c_1-o(1),c_2+o(1)]$.
\item Fix any positive $\d$. Then $\g^k\in[T_k-\d, T_k]$ provided $k$ is large enough.
\end{itemize}
\end{prop}
\begin{proof}
The first part of the statement follows immediately from the definition of $f$.

To show the second part, note that from Proposition~\ref{propfn} we have
$$
f_{n_k-1}(\g^k)<\frac 12=f_{n_k}(T_k)\le f_{n_k}(\g^k)\
\Longrightarrow  \ \g^k\le T_k.
$$
On the other hand, for any small positive $\d$
$$
f_{n_k}(T_k-\d)\ge f_{n_k}(T_k) e^{r c_1 \d}
=\frac{e^{r c_1 \d} }{2},
$$
hence from Proposition~\ref{propfn}
$$
f_{n_k-1}(T_k-\d)=f_{n_k}(T_k-\d)-\left[f_{n_k}(T_k-\d)-f_{n_k-1}(T_k-\d)\right]
\ge 
\frac{e^{r c_1 \d} }{2} - e^{-r c_1 (T_k-\d)}
\ge \frac 12 +\frac{r c_1 \d}{2} - e^{-r c_1 (T_k-\d)}.
$$
For $k$ large enough, $\g^k$ is also large, and so is $T_k\ge \g^k$, which yields that $ e^{-r c_1 (T_k-\d)}$ is very small, yielding that the RHS of the above expression is larger than $1/2$, so
$$
f_{n_k-1}(T_k-\d) >\frac12> f_{n_k-1}(\g^k)\ \Longrightarrow\ T_k-\d<\g^k.
$$
by monotonicity of $f_{n_k-1}(\cdot)$. Consequently, $\g^k\in[T_k-\d, T_k]$ for all large $k$.
\end{proof}

\begin{coro}\label{cor1}
For any $\eps>0$
\begin{align}\label{eqNnu}
\begin{split}
\P\left(N^\gg( (1-\eps)\, \g^k )\ge n_k\right)&\to 0,\\
\P\left(N^\gg( (1+\eps)\, \g^k )< n_k\right)&\to 0
\end{split}
\end{align}
as $k\to\infty$.
\end{coro}
\begin{proof}
Immediately  follows from Propositions~\ref{propbounds} and~\ref{propTnu}.
\end{proof}
Now we are ready to present the proof of our main result.

\begin{proof}[Proof of Theorem~\ref{thm_main}.]
Fix $k\ge 1$ and define the new {\em\textcolor{blue}{blue}} process  
$S^\bb$, which takes values in $\{0,1\}^{\Z_+}$, as follows. Let 
$\tau_{n_k}^{(1)}$ be the first time when the fire process reaches $n_k$. Up to this time $S^\bb$ completely coincides with the fire process $S^\ff$. However, at this hitting time, we set $S^\bb(x,\tau_{n_k}^{(1)})=0$ for all $x$ (not only just those reachable from $0$). After this, the blue process evolves again exactly like the fire process on the same probability space using the same collection of the Poisson process ${\cal P}_x$ (please see the beginning of the paper), until the next time when $n_k$ is burnt, at which point the blue process restarts again with zeroes everywhere. Thus $S^\bb(x,t)$ is a renewal process with the renewal times $\tau_{n_k}^{(i)}$, $i=1,2,\dots$, where $\tau_{n_k}^{(i)}$ are the consecutive times when the fire process reaches $n_k$. Moreover, the blue process and the fire process are trivially coupled such that
\begin{align*}
S^\bb(x,t)&=S^\ff(x,t)\quad\text{for }t<\tau_{n_k}^{(1)} \text{ and all }x;\\
S^\bb(x,t)&=S^\ff(x,t)\quad\text{for all }t
\text{ and }x\le n_k;\\
S^\bb(x,t)&\le S^\ff(x,t)\quad\text{for all }t
\text{ and }x> n_k.
\end{align*}
Let $\rho_i$ ($\rho_i^\ff$ respectively), $i=1,2,\dots$, be the right-most point burnt by the blue process (the fire process respectively) at time $\tau_{n_k}^{(i)}$. Define the {\em inter-arrival times} 
$\Delta_{n_k}^{(i)}=\tau_{n_k}^{(i)}-\tau_{n_k}^{(i-1)}$, $i=2,3,\dots$, with the convention that $\Delta_{n_k}^{(1)}=\tau_{n_k}^{(1)}=\tau_{n_k}$.
From the definition of the blue process, it follows that the random variables $\xi_i=\left(\Delta_{n_k}^{(i)},\rho_i\right)$, $i=2,3,\dots$ are i.i.d.

By construction $\rho_1=\rho_1^\ff$, that is the rightmost burnt point is the same for the fire and the blue processes when the fire process reaches $n_k$ for the first time; however, this does not necessarily hold for the times $\tau_{n_k}^{(i)}$, $i\ge 2$. At the same time we can claim
\begin{lemma}\label{cla2}
We have
\begin{itemize}
\item $\rho_i \le \rho_i^\ff$ for all $i$;
\item on the event $\rho_i\le \rho_{i-1}$ we have $\rho_i^\ff=\rho_i$.
\end{itemize}
\end{lemma}
\begin{proof}[Proof of Lemma]
The first part of the lemma trivially follows from the fact that $S^\bb(x,t)\le S^\ff(x,t)$ for all $x$ and $t$. To show the second statement, note that at time $s=\tau^{(i-1)}_{n_k}$ we have
$$
S^\ff\left(x,\tau_{n_k}^{(i-1)}\right)
=
S^\bb\left(x,\tau_{n_k}^{(i-1)}\right)
=0\quad\text{for }x=0,1,\dots,\rho_{i-1}.
$$
Until the time $\tau_{n_k}^{(i)}$ the fire and the blue processes coincide on $[0,\rho_{i-1}]$.
Since $\rho_{i}\le \rho_{i-1}$, 
$$
S^\bb\left(x,\tau_{n_k}^{(i)}-0\right)=0
\text{ for } x=\rho_i,\rho_i-1,\dots,\rho_i-(r-1)
$$
where $S^\bb(x,t-0)$ denotes the state of $x$ in the infinitesimal moment just before time $t$. However, since the blue and the fire process coincide on $[0,\rho_{i-1}]$, this also means that
$$
S^\ff\left(x,\tau_{n_k}^{(i)}-0\right)=0
\text{ for } x=\rho_i,\rho_i-1,\dots,\rho_i-(r-1)
$$
and thus the fire cannot reach any point beyond $\rho_i$.
\end{proof}

Fix now an $\eps>0$ so small that 
\begin{align}\label{eqgamma}
\g (1+\eps)<2\, (1-\eps).
\end{align}
For each $j=1,2,\dots$, let $S^\gg_j$ denote a version of the green process which is reset everywhere to zero at time $\tau_{n_k}^{(j-1)}$, i.e., we have $S^\gg_j\left(x,\tau_{n_k}^{(j-1)}\right)=0$ for all $x\in \Z_+$, with the convention $\tau_{n_k}^{(0)}\equiv 0$; thus
$S^\gg_j\equiv S^\gg_1$. Define 
\begin{align*}
E_{k,j}=\{&\text{for the green process $S^\gg_j$ there exists a subinterval of } (0, n_{k+1}]
\\
&\left. \text{  of length $r$ that has no Poisson arrivals during the time }\left[0,\Delta_{n_k}^{(j)}+\Delta_{n_k}^{(j+1)}\right]\right\}
\end{align*}
as well as $\a_k= \P\left(E_{k,j}\right)$; note this probability does not depend on $j$ since $\Delta_{n_k}^{(j)}$ are i.i.d.
\begin{lemma}\label{cla}
$\a_k\to 0$ as $k\to\infty$.
\end{lemma}
\begin{proof}[Proof of Lemma]
Recall $\tau_x^\gg$ from Definition~\ref{def1} and apply Corollary~\ref{cor1}. We have  shown
that
\begin{align}\label{eqtaubaunds}
\P(|\tau_{n_k}^\gg-\g^k|>\eps \g^k )\to 0
\text{ as }k\to\infty.
\end{align}
and  that
\begin{align}\label{eqtaunotsmall}
\P(\tau_{n_k}<\g^k(1-\eps))\le   \P\left(\tau^\gg_{n_k}<\g^k(1-\eps)\right)=o(1)
\end{align}
where $o(1)\to 0$ as $k\to\infty$. From~(\ref{eqgamma}), (\ref{eqtaubaunds}) and~(\ref{eqtaunotsmall})  it follows that 
\begin{align*}
\P(E_{k,j})& \le \P\left(E_{k,j} \cap \{\Delta_{n_k}^{(j)}+\Delta_{n_k}^{(j+1)}>2(1-\eps)\g^k\}\right)+2\,\P(\tau_{n_k}<(1-\eps)\g^k)
\\ & 
\le\P\left(\tau_{n_{k+1}}^\gg> 2 (1-\eps)\g^k\right)+o(1) 
\le\P\left(\tau_{n_{k+1}}^\gg> \g (1+\eps) \g^k \right)+o(1) 
\\
&=\P\left(\tau_{n_{k+1}}^\gg>  (1+\eps) \g^{k+1} \right)+o(1)
\le o(1)+o(1) \nonumber
\end{align*}
where $o(1)\to 0$ as $k\to\infty$.
\end{proof}
\begin{coro}\label{cla3}
For $j=1,2,\dots $ we have
$$
\left\{\rho_{j+1}^\ff\ge \rho_{j}^\ff \text{ and }\rho_{j+1}^\ff<n_{k+1}\right\}
\subseteq E_{k,j}.
$$
\end{coro}
\begin{proof}
The result follows from the construction of $S^\gg_j$ and $S^\ff$ which are both functionals of the collection of the Poisson processes ${\cal P}_x$, $x\in \Z_+$; in particular,
$$
S^\ff\left(x,t\right)\ge  S^\gg_j\left(x,t\right)
$$
for $t\in[\tau_{n_k}^{(j)},\tau_{n_k}^{(j+1)}]$ and $x\in [\rho_j^\ff,n_{k+1}]$.
\end{proof}

One of the fires that burns $n_k$ will eventually burn $n_{k+1}$ as well, so we can write
\begin{align*}
\tau_{n_{k+1}}&=\Delta^{(1)}_{n_k}
+\Delta^{(2)}_{n_k} 1_{\left\{\rho_1^\ff<n_{k+1}\right\}} 
+ \Delta^{(3)}_{n_k} 1_{\left\{\rho_1^\ff<n_{k+1},\rho_2^\ff<n_{k+1}\right\}}
+\dots
=\Delta^{(1)}_{n_{k}}+\sum_{i=1}^\infty \Delta^{(i+1)}_{n_{k}} 1_{A_i}
\end{align*}
where
$$
A_i=A_i^{(k)}=\bigcap_{j=1}^{i}\left\{\rho_j^\ff<n_{k+1}\right\}
$$
is a decreasing sequence of events. In other words, $A_i$ corresponds to the event that during the first~$i$ fires at point $n_{k}$, point $n_{k+1}$ has not yet been burnt. 
Since~$\Delta^{(i)}_{n_k}$, $i=1,2,\dots$, are i.i.d., $\Delta^{(1)}_{n_k}\equiv\tau_{n_k}$, and~$\Delta^{(i+1)}_{n_k}$ is independent of $A_i$ for each $i$, 
\begin{align}\label{eqEtaubis}
\E \left(\tau_{n_{k+1}}\right)=
\E \tau_{n_k}\cdot\left[1+\sum_{i=1}^\infty \P(A_i)\right].
\end{align}

Our task now is to estimate $\P(A_i)$ where $i\ge 1$.  
Let $y=(y_0,y_1,\dots,y_i)$ be a sequence of positive numbers of length $i+1$ with the convention that $y_0\equiv+\infty$. We say that $y_j$, $j=1,2,\dots,i-1$, is {\em a weak local minimum} if 
$$
H_j(y)=\{y_j\le  \min(y_{j-1},y_{j+1})\}
$$
holds.

For each such sequence $y$  there exists a non-negative integer $\nu=\nu(y)\ge 0$, which is a lower bound on the number of local minima, and the increasing sequence of indices $s(y)=(s_1(y),s_2(y),\dots,s(y)_{\nu(y)})$ with the property 
\begin{itemize}
\item $s_1(y)=\inf\{j\ge 1:\ H_j(y)\text{ occurs}\}$ is the index of the first weak local minimum in the sequence~$y$;
\item $s_{m+1}(y)$ is the index of the first local minimum in the sequence $y$ with the index of at least $s_m(y)+3$.
\item $s_{\nu}(y)+3>i-1$, or there is no local minimum with the index more than or equal to $s_{\nu}(y)+3$.
\end{itemize}
Set $\nu(y)=0$ if no weak local minima exist in $y$.
(To illustrate this concept, consider the  example where $ y=(\infty,3,2,4,1,3,2,5)$, then $\nu=2$ and $(s_1,s_2)=(2,6)$. Note that  the number of local minima here is $3>\nu$, since the middle local minimum is too close to the first one.)

Let now 
$$
y=(+\infty,\rho_1,\dots,\rho_i),
\qquad \sigma=s(y)
$$ 
where $\rho$s are defined in the beginning of the proof of the theorem, and observe that the second condition in the definition of $s_m(y)$ ensures that $\sigma_{m+1}\ge \sigma_m+3$, so the triples
$(\sigma_m-1,\sigma_m,\sigma_m+1)$ are non-overlapping for distinct $m$s, implying, in turn, that the triples
 $(\rho_{\sigma_m-1}$, $\rho_{\sigma_m}$, $\rho_{\sigma_m+1})$ are {\em independent} for distinct $m$s.
As a result,
$$
\P(\sigma_{m+1}=\sigma_m+3 \| \sigma_1,\dots,\sigma_m<i-3)\ge \P\left(\rho_{\sigma_m+3}\le \min(\rho_{\sigma_m+2},\rho_{\sigma_m+4})\right)\ge \frac 13.
$$
due to the symmetry between $\rho_{\sigma_m+2},\rho_{\sigma_m+3},\rho_{\sigma_m+4}$. Since $\rho_i$s are i.i.d., we conclude therefore that $\nu(y)$ is stochastically larger than a ${\rm Binomial}(\lfloor (i+1)/3\rfloor,1/3)$ random variable.

By Lemma~\ref{cla2} $\rho_{j+1}^\ff\ge \rho_{j+1}$
and moreover on the event $H_j(\rho)$ we have $\rho_j^\ff=\rho_j$ since  $\rho_{j}\le \rho_{j-1}$. Since also $(\rho_{j}^\ff=)\rho_{j}\le \rho_{j+1} (\le \rho^\ff_{j+1})$ on $H_j(\rho)$, by Lemma~\ref{cla2} and Corollary~\ref{cla3}, 
$$
\P\left(\rho_{j+1}^\ff < n_{k+1},\ H_j\right)
\le
\P\left(\rho_j^\ff \le \rho_{j+1}^\ff < n_{k+1},\ H_j\right)
\le
\P(E_{k,j},\ H_j)\le\P(E_{k,j})\le \a_k\to 0.
$$

By the law of iterated expectations
$$
\P(A_i,\ \nu(y)=m)
=\int_{y:\ \nu(y)=m}\P(A_i\| (\rho)_i=y){\rm d}\P(\rho=y) 
\le \sup_{y:\ \nu(y)=m} \P(A_i\| (\rho)_i=y)
$$
where $(\rho)_i:=(+\infty,\rho_1,\rho_2,\dots,\rho_i)$, and
the supremum is taken over all sequences $y$ of length $i+1$ with $y_0=+\infty$ and all other elements being positive. However, when $\nu(y)=m$ for any realization of $s(y)=(\sigma_1,\dots,\sigma_m)$ we have 
\begin{align}\label{eqakm}
\P(A_i \| (\rho)_i=y)
&=\P\left( \bigcap_{j=1}^i \left\{\rho_{i}^\ff < n_{k+1}\right\}
 \| (\rho)_i=y \right)
 \nonumber
\\ &
\le 
\P\left( \bigcap_{l=1}^m \left[
\left\{\rho_{\sigma_l+1}^\ff < n_{k+1}\right\}
\cap H_{\sigma_l}(\rho)
\right]
 \|(\rho)_i=y \right)
\le \a_k^m
\end{align}
due to the independence of the triples $\{\rho_{l-1},\rho_l,\rho_{l+1}\}$ for different $l\in s(y)$. Consequently,
\begin{align}\label{eqAii+mbis}
\P(A_i)&=\P\left(A_i \cap \left(\bigcup_{m=0}^\infty\{\nu(y)=m\}\right)\right)
=\sum_{m=0}^\infty \P(A_i,\nu(y) =m)\nonumber 
\\
&
\le \P(\nu(y)=0) 
+\sum_{m=1}^{\lfloor i/10\rfloor+1}  \P(A_i,\nu(y)=m)
+ 
\sum_{m=\lfloor i/10\rfloor+2}^{\infty}  \P(A_i,\nu(y)=m)
\nonumber
\\
&
 \overset{  \text{by~\eqref{eqakm}} }{\le} 
\P(\nu(y)=0)+\a_k\,\P\left(1\le \nu(y)\le \frac{i}{10}+1\right) +\sum_{m=\lfloor i/10\rfloor+2}^{\infty}  \a_k^m
\nonumber
\\
&
\le \P(\nu(y)=0)+  \a_k\,2 e^{-c_4 i} +\frac{\a_k^{1+i/10}}{1-\a_k}
\end{align}
since
$$
\P\left( \nu(y)\le \frac{i}{10} +1\right) 
\le \P\left({\rm Binomial}\left(\left\lfloor \frac{i+1}3\right\rfloor,\frac13\right) \le \frac{i}{10} +1\right)\le   2\,e^{-c_4 i}.
$$
for some $c_4>0$ by the large deviation theory (see, e.g.~\cite{denH}).
Finally, observe that on $i\ge 2$
$$
\P(\nu(y)=0)=\P(\rho_1>\rho_2>\dots >\rho_i)\le \frac{1}{i!}
$$
due to the symmetry of all permutations of $[1,2,\dots,i]$.

Trivially, $\P(A_1)\le 1$, so by~\eqref{eqAii+mbis}
$$
\sum_{i=1}^\infty \P(A_i)\le  1+\sum_{i=2}^\infty 
\left[\frac 1{i!}+\a_k\left(2 e^{-c_4 i} +\frac{\a_k^{i/10}}{1-\a_k}\right) \right]=e-1+o(1).
$$
where $o(1)\to 0$ as $k\to\infty$, since $\a_k\to 0$.
Therefore~(\ref{eqEtaubis}) gives
\begin{align}\label{eq:Cbis}
\E \left(\tau_{n_{k+1}}\right)&=
\E \tau_{n_k}\cdot\left[1+\sum_{i=1}^\infty \P(A_i)\right]
\le 
(e+o(1))\cdot \E \tau_{n_k},
\end{align}
yielding that for any fixed $\d>0$
$$
\E \left(\tau_{n_k}\right) \le   (e+\d)^k
$$
for all sufficiently large $k$.

For each $x>0$ one can find a unique $k$ such that $n_{k-1}< x\le  n_{k}$, and by Proposition~\ref{propTnu},
$$
k= \log_\g\left(\frac{\log x }{r\tilde c_k}\right)
+O(1)= \log_\g\left(\log x \right)+O(1),
$$
whence
$$
\E \tau_x \le \E \tau_{n_k}
\le 
\left(\log x\right)^{\log_\g (e+2\d)}
$$ 
where the power can be made arbitrary close to $(\log 2)^{-1}=1.442695\dots $ by choosing $\g\uparrow 2$ and $\d\downarrow 0$. 
This proves Theorem~\ref{thm_main}.
\end{proof}

\section{Continuous tree model on $\R_+$}
\label{sec-cont}
We want to get some estimates for the ``green'' process in case of the continuous space model, i.e.\ the one where by time $t\ge 0$ we have a Poisson point process on~$\R_+$ (=the set of occupied sites) in space--time with intensity ${\rm d}x\otimes {\rm d}t$, and each point is a center of a circle of radius $1$.

We say  that two sites $x$ and $y$ of the Poisson process are {\it connected} if $|x-y|\le 1$. We assume that $0$ is always occupied. With these definitions we can also define the cluster of occupied sites containing zero, which will be the subset of points of the Poisson process $x_1=x_1(t)$, $x_2=x_2(t)$, $\dots$, $x_{n(t)}=x_{n(t)}(t)$ such that
$$
x_1\le 1,\ x_2-x_1\le 1, \ \dots, \ x_n-x_{n-1}\le 1, \ x_{n+1}-x_n> 1 
$$
where $n=n(t)$ is the number of sites in the cluster. Let also $N^\ff(t)=x_{n(t)}(t)$ be the location of the right-most site in the cluster, and $\tau_x^\ff $, $x>0$, be the smallest positive time for which $x\le N^\ff(t)$; thus 
$$
\{\tau_x^\ff >t\}=\{N^\ff(t)<x\}.
$$
We will find the estimate for $N^\gg(t)$ and $\tau_x^\gg$; some similar results can be found in literature, see e.g.\ Proposition~5.2 in~\cite{ASMU}.

\begin{prop}\label{propcont}
For the green process in the continuous model on $\R_+$ we have
\begin{align*}
\P\left(N^\gg(\log x)\ge x\right)&=o(1),\\
\P\left(N^\gg(\log x +3\log \log x\right)\le  x )&=o(1).
\end{align*}
\end{prop}
\begin{proof}
Let $W_1=x_1$, $W_i=x_{i}-x_{i-1}$, $i\ge 2$. Then~$W_i$ are i.i.d.\ exponentially distributed random variables with rate~$t$. We can compute the Laplace transform of $N^\gg(t)$ as follows:
\begin{align}\label{eq:laplaceGreenc}
\E\left[ e^{-\l N^\gg(t)} \right] &= \E\left[ e^{-\l \sum_{i=1}^{n(t)} W_i} \right] 
= \sum_{m=0}^{\infty} \E\left[ e^{-\l \sum_{i=1}^{m} W_i} {\bf 1}_{\{n(t)=m\}} \right] \nn
\\ \nn
&= \sum_{m=0}^{\infty} \E\left[ \prod_{i=1}^m e^{-\l W_i}  {\bf 1}_{\{ W_i\leq 1\}}  
\times  {\bf 1}_{\{ W_{m+1}> 1\}}  \right] 
= \sum_{m=0}^{\infty} \E\left[  e^{-\l w_1}  {\bf 1}_{\{ w_1\leq 1\}}  \right]^m
\times  \P\left[ w_1>1\right] 
\\ 
&= \sum_{m=0}^{\infty} \left[  \frac{t}{t+\l} \left(1-e^{-\l-t}\right) \right]^m \times  e^{-t} =  \frac{(\l+t)e^{-t}}{\l+te^{-\l-t}}.
\end{align}
Using the Taylor expansion at $\l =0$ we get the first moments of $N^\gg(t)$,
$$
\E N^\gg(t)= \frac{e^t-1-t}{t} \;,\qquad {\rm Var}(N^\gg(t))= \frac{e^{2t}-1-2te^t}{t^2}=
\frac{(e^t-1-t)^2+2(e^t-1-t-t^2/2)}{t^2}.
$$
It is easy to check that 
$$
\lim_{t \to \infty} \E\, 
e^{- \displaystyle\frac{\l N^\gg(t)}{\E N^\gg(t)}}   = \frac1{1+\l}
$$ 
for $\Re e \l>-1$,
thus 
$$
te^{-t} N^\gg(t) \longrightarrow {\rm Exponential\; mean\ }  1, \; {\rm in\ law.}
$$
Consequently, as $t\to\infty$,
$$
\P(N^\gg(t)>e^t)=\P(te^{-t} N^\gg(t) >t)=o(1)\quad \text{and}
\quad
\P\left(N^\gg(t)<e^t/t^2\right)=\P\left(te^{-t} N^\gg(t) <1/t\right)=o(1).
$$
This, in turn, implies
$$
\P(N^\gg(\log x)>x)=o(1)
$$
and
$$
\P\left(N^\gg(\log x+3\log\log x)<x\right)\le
\P\left(N^\gg(\log x+3\log\log x)< \frac{x\log^3 x}{(\log x+3\log\log x)^2}\right)=o(1)
$$
implying the statement of the proposition.
\end{proof}

The following statement can be proven following verbatim the lines of the proof in Section~\ref{Sec_couplingfire} with $x_k=\g^k$, since the estimate~\eqref{eqNnu} is ensured by Proposition~\ref{propcont}.
\begin{thm}\label{thm_cont}
For the  continuous model of forest fire  on $\R_+$ we have
that for any $\delta>0$ 
$$
\E \tau_x\le  (\log x)^{(\log 2)^{-1}+\delta}.
$$ 
for all $x$ large enough.
\end{thm}

%

\end{document}